\documentclass[12pt]{amsart}
\usepackage{amssymb,amsmath,amsthm,amsfonts,amsopn,url,bbm,amscd,enumerate}
\theoremstyle{plain}
\newtheorem{thm}{Theorem}
\newtheorem{prop}{Proposition}
\newtheorem{lemma}{Lemma}
\newtheorem{cor}{Corollary}
\theoremstyle{definition}
\newtheorem{definition}{Definition}
\newcommand{\field}[1]{\mathbbm{#1}}

\newcommand{\Z}{\field{Z}}

\newcommand{\ideal}[1]{\mathfrak{#1}}
\newcommand{\m}{\ideal{m}}
\newcommand{\n}{\ideal{n}}
\newcommand{\p}{\ideal{p}}

\newcommand{\aaa}{\ideal{a}}
\newcommand{\ffunc}[1]{\mathrm{#1}}
\newcommand{\func}[1]{\mathop{\rm #1}}

\newcommand{\bx}[1]{{\mathbf{x}^{[#1]}}}

\newcommand{\bfx}{\mathbf{x}}

\newcommand{\HKmult}{\ffunc{e_{\rm H K}}}
\newcommand{\HKmul}[1]{\mathrm{e}_{\rm H K}^{#1}}
\newcommand{\sptc}[1]{{#1}^{*\rm sp}}
\newcommand{\z}{\mathbf{z}}
\begin{document}
\title{A length characterization of $*$-spread}
\author{Neil Epstein and Adela Vraciu}
\address{Department of Mathematics, University of Michigan, Ann Arbor, MI 48109}
\address{Department of Mathematics, University of South Carolina, Columbia SC 29205}
\email{neilme@umich.edu; vraciu@math.sc.edu}
\date{\today}
\begin{abstract} 
The $*$-spread of an ideal is defined as the minimal number of generators of an ideal which is minimal with respect to having the same tight closure as the original ideal. We prove an asymptotic length formula for the $*$-spread.               
\end{abstract}
\maketitle
\section{Introduction}
Several closure operations for ideals in a commutative Noetherian rings have been studied by numerous authors; among those closures, we mention integral closure, tight closure, Frobenius closure, and plus closure.

For each of the above-mentioned closure operations, a corresponding notion of spread can be defined as the minimal number of generators of a minimal reduction with respect to that operation. The fact that the minimal number of generators is independent of the choice of the reduction is well-known in the case of the integral closure (\cite{NR}), easy to see in the case of Frobenius closure, and recently proved (\cite{Epst}) in the case of tight closure.

Note that in most cases, these spreads can be characterized
asymptotically in terms of length, and without reference to corresponding reductions of the ideal.
  In the case of
integral closure, bar-spread is equal to analytic spread (provided the residue field is infinite): \[
\ell(I) = \ell^-(I) = \deg_n \dim_k \lambda(I^n / \m I^n) + 1 = \deg_n \dim_k \mu(I^n) + 1.
  \footnote{Here, $\lambda$ stands for length, and $\mu$ is the minimal-number-of-generators function.}\] The $F$-spread of $I$ is the eventual minimal number of generators
of high Frobenius powers of $I$.  That is, \[
\ell^F(I) = \lim_{q \rightarrow \infty} \lambda(I^{[q]} / \m I^{[q]}) = \lim_{q \rightarrow \infty} \mu(I^{[q]}).
\]  Finally, the $+$-spread of an ideal $I$ in a henselian local domain $R$ is the eventual minimal 
number of generators of $I$ expanded to domains which are integral extensions of $R$: \[
\ell^+(I) = \func{\rm colim}_{\stackrel{(R,\m) \subseteq (S,\n)}{\text{ integral ext. domain}}} \lambda(I S / \n I S) = 
\func{\rm colim}_{\stackrel{(R,\m) \subseteq (S,\n)}{\text{ integral ext.  domain}}} \mu(I S).
\]
The main result of this paper is an asymptotic characterization of $*$-spread (the spread corresponding to tight closure) in terms of length
\footnote{If one assumes that $J$ is $\m$-primary and $*$-independent, then the theorem is proved in \cite[Theorem 3.5(a)]{Vr}}:
\begin{thm}\label{thm:*spreadlength}
Let $(R,\m,k)$ be an analytically irreducible excellent local ring of characteristic $p>0$ and Krull dimension $d$ such that
$k = \kappa(\bar{R})$.  Let $J$ be a proper ideal,
and let $\aaa$ be an $\m$-primary ideal.  Then for $q_0 \gg 0$, 
\begin{equation}\label{eq:*spreadformula}
\ell^*(J) = \frac{1}{\HKmult(\aaa)} \lim_{q \rightarrow \infty}
 \frac{\lambda(J^{[q q_0]} / \aaa^{[q]} J^{[q q_0]})} {q^d}.
\end{equation}
In particular, if $\lambda(R/J) < \infty$, then \begin{equation}\label{eq:*spreadHKformula}
\ell^*(J) = \frac{\HKmult(\aaa J^{[q_0]}) - \HKmult(J^{[q_0]})}{\HKmult(\aaa)}
\end{equation}
\end{thm}

Here, $\kappa(\bar{R})$ stands for the residue field of the normalization $\bar{R}$ of $R$ (which is a local domain, due to the analytic irreducibility of $R$; a proof of this fact can be found in \cite[Lemma 4.3]{Epst}, although it has been known as folklore before).

As an application, we get a result which connects the rationality of the Hilbert-Kunz multiplicity for the ideals $I, J$, and $IJ^{[q]}$, where $I$ and $J$ are $\m$-primary ideals (Proposition~\ref{prop:HKrat}). We also prove a change of base formula for $*$-spread under flat local homomorphisms (Proposition~\ref{prop:*spreadBC}).
\section{Preliminaries}
We review some of the notions and results that are used in the proof of our main result. Note that in this paper $p$ always stands for the characteristic of $R$, and $q, q', q_0, q_1, q_2$ stand for various powers of $p$.

\begin{definition}
Let $R$ be a Noetherian local ring of characteristic $p>0$, let $f_1, \dotsc, f_\ell \in R$. We say that $f_1, \dotsc, f_\ell$ are {\it $*$-independent} if $f_i \notin (f_1, \dotsc, \hat{f_i}, \dotsc f_\ell)^*$ for all $i\in [1,\ell]
\cap \Z$.

We say that an ideal $I\subset R$ is $*$-independent if can be generated by $*$-independent elements. 
If $R$ is local, excellent, and analytically irreducible, this is equivalent to every minimal system of generators being $*$-independent \cite[Proposition 3.3]{Vr}.
\end{definition}

\begin{definition}
Let $R$ be a Noetherian local ring of characteristic $p>0$, $I, K\subset R$ ideals. We say that $K$ is a {\it $*$-reduction of $I$} if $K \subseteq I\subseteq K^*$. We say that $K$ is a {\it minimal $*$-reduction of $I$} if it is minimal with this property.
\end{definition}

 Note that this is equivalent to being a  $*$-reduction such that every minimal set of generators is $*$-independent \cite[Proposition 2.1 and Lemma 2.3]{Epst}. 
 Minimal $*$-reductions exist by \cite[Propositions 2.1 and 2.3]{Epst}.

\begin{definition}
Let $(R, \m, k)$ be an excellent analytically irreducible local domain of characteristic $p>0$, $I\subset R$ an ideal. The {\it $*$-spread} of $I$, denoted $\ell^*(I)$, is the minimal number of generators of a minimal $*$-reduction of $I$. The fact that this number is independent of the choice of a minimal $*$-reduction is Theorem 5.1 in \cite{Epst}. 
\end{definition}

\begin{definition}\cite{Vr}
Let $(R, \m)$ be a Noetherian local ring of characteristic $p>0$, $x \in R$, $I\subset R$ an ideal. We say that $x \in I^{*sp}$, the {\it special tight closure of $I$}, if there exists $q_0 =p^{e_0}$ such that $x^{q_0} \in (\m I^{[q_0]})^*$. 

Note that one can replace $\m I^{[q_0]}$ by $\aaa I^{[q_0]}$ in the above definition for any $\m$-primary ideal $\aaa$, by suitably increasing $q_0$. 
\end{definition}

The following result was proved in \cite{Epst}:
\begin{thm}\label{thm: specialtc}
Let $(R, \m, k)$ be an excellent analytically irreducible local domain of characteristic $p>0$. Assume that $k=\kappa(\bar{R}$. Then for any proper ideal $I$ of $R$, there exists $q'$ such that 
$$
(I^*)^{[q']} \subseteq I^{[q']} +\sptc{(I^{[q']})}
$$
\end{thm}

We will also use the following result of \cite{Ab}, which we will refer to as the {\it Colon Criterion}:  
\begin{prop}\label{prop: coloncrit}
Let $(R, \m)$ be an excellent analytically irreducible local domain, let $I\subset R$, and $x \notin I^*$. Then there exists $q_0$ such that \newline
$I^{[q]}: x^q \subset \m^{[q/q_0]}$ for all $q \ge q_0$.
\end{prop}

\begin{definition}
Let $(R, \m)$ be a Noetherian local ring of characteristic $p>0$, and $I\subset R$ an $\m$-primary ideal. The {\it Hilbert-Kunz multiplicity} of $I$ is
$$
e_{HK}(I)=\lim_{q \rightarrow \infty} \frac{\lambda(R/I^{[q]})}{q^d},
$$
where $d$ is the Krull dimension of $R$.
\end{definition}
The Hilbert-Kunz multiplicity was identified as the coefficient of the leading term of the Hilbert-Kunz function $\lambda(R/I^{[q]})$ in \cite{Mo}, and it turned out to be an important tool in the study of tight closure, due to \cite[Theorem 8.17]{HH}, which asserts that two $\m$-primary ideals
$I\subseteq J$ have the same tight closure if and only if they have the same Hilbert-Kunz multiplicity.

\section{Proof of the Main Result}

Before we prove our main result, Theorem~\ref{thm:*spreadlength}, we need some preliminary results.

\begin{lemma}\label{lem:*indparam}
Let $R$ be a Noetherian local ring of characteristic $p>0$, let $f_1, \dotsc, f_\ell$
be $*$-independent elements generating an ideal $K$, and
let $\bfx = x_1, \dotsc, x_n$ be parameters modulo $K$.
Then there is some positive integer
$t$ such that $f_1, \dotsc, f_\ell, x_1^t, \dotsc, x_n^t$ are $*$-independent.
\end{lemma}

\begin{proof}
First note that $x_j^t \notin (f_1, \dotsc, f_\ell, x_1^t, \dotsc, \hat{x_j^t}, \dotsc, x_n^t)^*$ for $1 \leq j \leq n$, $\forall \, t$
since the $x_i$ are parameters mod $K$, so the heights of the latter ideal and $K + (\bfx^t)$ do not match modulo $K$.

Now pick some $1 \leq i \leq \ell$, and suppose $f_i \in (f_1, \dotsc, \hat{f_i}, 
\dotsc, f_\ell, x_1^t, \dotsc, x_n^t)^*$ for all $t$.  Then since each of these ideals contains the next, \[
f_i \in \bigcap_{t \geq 1} (f_1, \dotsc, \hat{f_i}, \dotsc, f_\ell, x_1^t, \dotsc, x_n^t)^* 
= (f_1, \dotsc, \hat{f_i}, \dotsc,  f_\ell)^*,
\]
which contradicts the $*$-independence of the $f_j$.  (The equality holds essentially because of the Krull intersection theorem.)
Thus for each $i$ with $1 \leq i \leq \ell$,
there exists an integer $t_i$ with \[
f_i \notin (f_1, \dotsc, \hat{f_i}, \dotsc, f_\ell, x_1^{t_i}, \dotsc, x_n^{t_i})^*.
\]  Let
$t = \max_i \{ t_i \}$.  Then $K + (\bx{t})$ is a $*$-independent ideal.
\end{proof}

\begin{lemma}\label{lem:nottc}
Let $(R, \m)$ be a Noetherian local ring of characteristic $p>0$. Assume that $R$ has a weak test element $c$. Let $f_1, \dotsc, f_\ell \in R$ be $*$-independent, and $g_1, \dotsc, g_r \in \sptc{(f_1, \dotsc, f_\ell)}$.

Then $f_i \notin (f_1, \dotsc, \hat{f_i}, \dotsc, f_\ell, g_1, \dotsc, g_r)^*$ for all $i \in [1, \ell]$.
\end{lemma}
\begin{proof}
Let $c$ be a weak test element, and choose $q_0$ such that 
\begin{equation}\label{eq: sptc}
c g_j^{q q_0} = \sum_{k=1}^\ell a_{k j} f_k^{q q_0},
\end{equation}
with $a_{k j} \in \m^{[q]}$, and 
$$(f_1, \dotsc, \hat{f_i}, \dotsc, f_\ell)^{[qq_0]} : f_i^{qq_0}\subseteq \m^{[q]}.$$
The first equation holds by the definition of special tight closure, and the second by the colon criterion.

Assume by contradiction that $f_i \in (f_1, \dotsc, \hat{f_i}, \dotsc, f_\ell, g_1, \dotsc, g_r)^*$, so that
$$
c f_i^{q q_0}=\sum_{j=1}^r d_j g_j^{q q_0} + \sum_{\stackrel{1 \leq k \leq \ell}{k \neq i}} e_k f_k^{q q_0}.
$$
Multiplying this by $c$ and combining with equation~\ref{eq: sptc} yields
$$(c^2 - \sum_{j=1}^r d_j a_{i j}) \in (f_1, \dotsc, \hat{f_i}, \dotsc, f_\ell, z)^{[q q_0]} : f_i^{q q_0}
 \subseteq \m^{[q]}.$$
Since each $a_{i j} \in \m^{[q]}$, this implies that $c^2 \in \m^{[q]}$, which 
gives a contradiction since $q$ was taken to be arbitrarily large, but $c \in R^o$ was fixed.
\end{proof}

\begin{lemma}\label{lem:lengthinduction}
Let $(R,\m,k)$ be an excellent analytically irreducible local ring of characteristic $p>0$ such 
that $k = \kappa(\bar{R})$.  Let $I$ be a proper ideal which is not $\m$-primary,
let $\aaa$ be an $\m$-primary ideal, let $L$ be a minimal $*$-reduction of $I$, 
and let $z$ be a parameter modulo $I$ such that $L + (z)$ is a $*$-independent ideal.
Then there is some power $q_0$ of $p$ such that, \[
\lim_{q \rightarrow \infty} \frac{\lambda\left((I,z)^{[q q_0]} / \aaa^{[q]} (I,z)^{[q q_0]}\right)} {q^d}
 = \HKmult(\aaa) + \lim_{q \rightarrow \infty} \frac{\lambda\left(I^{[q q_0]} / \aaa^{[q]} I^{[q q_0]} \right)} {q^d}
\]
\end{lemma}

\begin{proof}
Let $L=(f_1, \dotsc, f_\ell)$ be a minimal generating
set of $L$, and let $I=(f_1, \dotsc, f_\ell, g_1, \dotsc, g_r)$ be a minimal generating set of $I$.
Such a minimal generating set exists by \cite[Lemma 2.2]{Epst}.
By Lemma~\ref{lem:*indparam}, we can pick $z$ such that $f_1, \dotsc, f_\ell, z$
are $*$-independent.  Thus there is a power $q_0$ of $p$ such that \begin{equation}\label{eq:fz*ind}
(L^{[q q_0]} : z^{q q_0}) + \sum_{i=1}^\ell (f_1, \dotsc, \hat{f_i}, \dotsc, f_\ell, z)^{[q q_0]} : f_i^{q q_0} \subseteq \aaa^{[q]}
\end{equation}
for all $q$, by the Colon criterion (Proposition~\ref{prop: coloncrit}).

Now, consider the following short exact sequence: \begin{equation}\label{eq:ses1}
0 \rightarrow \frac{\aaa^{[q]} I^{[q q_0]} + (z^{q q_0}) } {\aaa^{[q]} (I,z)^{[q q_0]}}
 \rightarrow \frac{(I,z)^{[q q_0]}} {\aaa^{[q]} (I,z)^{[q q_0]}}
 \rightarrow \frac{(I,z)^{[q q_0]}} {\aaa^{[q]} I^{[q q_0]} + (z^{q q_0})} \rightarrow 0.
\end{equation}

The first term is isomorphic
to $R / ((\aaa^{[q]} (I,z)^{[q q_0]}) : z^{q q_0})$.
Now, by the Colon Criterion, there is some $q_1$ such that $I^{[q]} : z^q \subseteq \m^{[q / q_1]}$
for all $q \geq q_1$.  Also, since $\aaa$ is $\m$-primary, there is some $q_2$ such that $\m^{[q_2]}
\subseteq \aaa$.  Let $q_0 \geq q_1 q_2$.

Let $u \in \left(\aaa^{[q]} (I,z)^{[q q_0]}\right) : z^{q q_0}$.  Then there is some $a \in \aaa^{[q]}$
such that \[
u - a \in (\aaa^{[q]} I^{[q q_0]}) : z^{q q_0} \subseteq I^{[q q_0]} : z^{q q_0} \subseteq
\m^{[q q_0 / q_1]} \subseteq \aaa^{[q q_0 / (q_1 q_2)]} \subseteq \aaa^{[q]}.
\]
Hence $u \in \aaa^{[q]}$.  The reverse containment is obvious, so \[
\aaa^{[q]} = \left(\aaa^{[q]} (I,z)^{[q q_0]}\right) : z^{q q_0}.
\]
Hence, the first term of (\ref{eq:ses1}) has length $\lambda(R/\aaa^{[q]})$.

For the third term of the sequence, we have: \begin{equation}\label{eq:blackbox}
\frac{(I,z)^{[q q_0]}} {\aaa^{[q]} I^{[q q_0]} + (z^{q q_0})}
 \cong \frac{I^{[q q_0]}} {I^{[q q_0]} \cap (\aaa^{[q]} I^{[q q_0]} + (z^{q q_0}))}
\end{equation}

\noindent \textbf{Claim:} We can choose $q_0$ to be large enough so that for any $q$, \[
\lim_{q \rightarrow \infty} \frac{\lambda\left(\left(I^{[q q_0]} \cap (\aaa^{[q]} I^{[q q_0]} + (z^{q q_0}))\right)
 / \aaa^{[q]} I^{[q q_0]}\right)} {q^d} = 0.
\]

\begin{proof}[Proof of Claim]
First note that the numerator of the above quotient of ideals equals
$\aaa^{[q]} I^{[q q_0]} + I^{[q q_0]} \cap (z^{q q_0})$.  Next, by Theorem~\ref{thm: specialtc}, there is some
$q_1$ such that \[
I^{[q_1]} \subseteq (L^*)^{[q_1]} \subseteq L^{[q_1]} + \sptc{(L^{[q_1]})}
\]
Hence, by replacing the $f_i$'s, the $g_j$'s, and $z$ by their $q_1$ powers, we may
assume that $I \subseteq L + \sptc{L}$.

After this replacement, then, there exist $h_i \in L$ and $g_i' \in \sptc{L}$
such that $g_i = g_i' + h_i$ for $1 \leq i \leq r$. We may replace the $g_i$ with the $g_i'$ and assume
without loss of generality that $g_i \in \sptc{L}$ for $1 \leq i \leq r$. By increasing $q_0$ if necessary, we may assume
$g_i^{q_0}\in (\aaa L^{[q_0]})^*$.

By Lemma~\ref{lem:nottc} we have $f_i\notin (f_1, \dotsc, \hat{f_i}, \dotsc, f_\ell, g_1, \dotsc, g_r, z)^*$.

Let $H_j := (g_1, \dotsc, g_j)$, where $1 \leq j \leq r$ (so $I = L + H_r$), and $H_0 := (0)$.  We show that 
$I^{[q q_0]} \cap (z^{q q_0}) \subseteq \aaa^{[q]} I^{[q q_0]} + H_r^{[q q_0]}$. Let $x \in I^{[q q_0]} \cap (z^{q q_0})$.

Then \[
x = t z^{q q_0} = \sum_{i=1}^\ell u_i f_i^{q q_0} + \sum_{j=1}^r v_j g_j^{q q_0},
\]
so that for each $1 \leq i \leq \ell$, \[
u_i \in (f_1, \dotsc, \hat{f_i}, \dotsc, f_\ell, g_1, \dotsc, g_r,z)^{[q q_0]} : f_i^{q q_0} \subseteq \aaa^{[q]}
\]
by the colon criterion (increasing $q_0$ if necessary).  Thus, $x \in \aaa^{[q]} I^{[q q_0]} + H_r^{[q q_0]}$, as claimed.

Note that \[
\lambda\left(\frac{I^{[q q_0]} \cap (\aaa^{[q]} I^{[q q_0]} + (z^{q q_0}))}{\aaa^{[q]} I^{[q q_0]}} \right)
 \leq \lambda\left(\frac{\aaa^{[q]} I^{[q q_0]} + H_r^{[q q_0]}} {\aaa^{[q]} I^{[q q_0]}} \right)
\]
by what we have shown immediately above.  Let $c$ be a test element. We have: 
$$
\lambda\left(\frac{\aaa^{[q]} I^{[q q_0]} + H_r^{[q q_0]}} {\aaa^{[q]} I^{[q q_0]}} \right)
 = \sum_{j=1}^r \lambda\left(\frac{\aaa^{[q]} I^{[q q_0]} + H_j^{[q q_0]}} {\aaa^{[q]} I^{[q q_0]} + H_{j-1}^{[q q_0]}} \right)
 = \sum_{j=1}^r \lambda\left(\frac{R}{(\aaa^{[q]} I^{[q q_0]} + H_{j-1}^{[q q_0]}) : g_j^{q q_0}}\right) 
$$
$$
 \leq \sum_{j=1}^r \lambda\left(\frac{R}{\aaa^{[q]} + (c)}\right) = r \lambda\left(\frac{R}{\aaa^{[q]} + (c)}\right)
= r \lambda_{(R/c)}\left(\frac{R / c} {((\aaa + (c)) / (c))^{[q]}}\right)
$$
The inequality is true because
$c g_j^{q q_0} \in \aaa^{[q]} L^{[q q_0]} \subseteq \aaa^{[q]} I^{[q q_0]}$.  Thus $\aaa^{[q]} + (c)
\subseteq \aaa^{[q]} I^{[q q_0]} : g_j^{q q_0}$, which proves the inequality.

The last term is $r$ times a Hilbert-Kunz function
over the $d-1$ dimensional ring $R/c$, 
hence bounded by
a constant times $q^{d-1}$, which proves the claim.
\end{proof}

At this point, taking limits of lengths over $q^d$ as $q \rightarrow \infty$ in (\ref{eq:blackbox}) 
gives: \[
\lim_{q \rightarrow \infty} \frac{\lambda(I^{[q q_0]} / (I^{[q q_0]} \cap (\aaa^{[q]} I^{[q q_0]} 
 \cap (z^{q q_0}) )))} {q^d} 
 = \lim_{q \rightarrow \infty} \frac{\lambda(I^{[q q_0]} / \aaa^{[q]} I^{[q q_0]})} {q^d},
\]
so that the exact sequence (\ref{eq:ses1}) yields that \begin{align*}
\lim_{q \rightarrow \infty} \frac{\lambda((I,z)^{[q q_0]} / \aaa^{[q]} (I,z)^{[q q_0]})} {q^d}
 &= \lim_{q \rightarrow \infty} \frac{\lambda(R / \aaa^{[q]}) +
 \lambda(I^{[q q_0]} / \aaa^{[q]} I^{[q q_0]})} {q^d} \\
 &= \HKmult(\aaa) + \lim_{q \rightarrow \infty} \frac{\lambda(I^{[q q_0]} / \aaa^{[q]} I^{[q q_0]})} {q^d}.
\end{align*}
\end{proof}

Now we begin the proof of Theorem~\ref{thm:*spreadlength}.

\begin{proof}
First suppose that $\lambda(R/J) < \infty$.

Let $K$ be a minimal $*$-reduction of $J$.  Consider the short exact sequences:
\begin{equation}\label{eq:firstses}
0 \rightarrow \frac{K^{[q q_0]}}{\aaa^{[q]} K^{[q q_0]}}
 \rightarrow \frac{J^{[q q_0]}}{\aaa^{[q]} K^{[q q_0]}}
 \rightarrow \frac{J^{[q q_0]}}{K^{[q q_0]}} \rightarrow 0
\end{equation}
and \begin{equation}\label{eq:secondses}
0 \rightarrow \frac{(\aaa J^{[q_0]})^{[q]}} {(\aaa K^{[q_0]})^{[q]}}
 \rightarrow \frac{J^{[q q_0]}} {\aaa^{[q]} K^{[q q_0]}}
 \rightarrow \frac{J^{[q q_0]}} {\aaa^{[q]} J^{[q q_0]}} \rightarrow 0.
\end{equation}

Since $J$ (and hence also $K$, since ideals with the same tight closure have
the same radical) is $\m$-primary, the length of the third term in (\ref{eq:firstses})
is the difference of the Hilbert-Kunz functions of $J$ and $K$.
Since these two have the same H-K multiplicity (since they have the same
tight closure), the limit as $q \rightarrow \infty$ of this difference
divided by $q^d$ is 0.  Hence the first and second terms are ``equal
in the limit''.

The same comment applies to the first term of the second short exact sequence,
since we have \[
\aaa J^{[q_0]} \subseteq \aaa (K^*)^{[q_0]} \subseteq \left(\aaa K^{[q_0]}\right)^*.
\]
Thus, the second and third terms of the second short exact sequence are
also ``equal in the limit''.  Hence by transitivity, \[
\lim_{q \rightarrow \infty} \frac{\lambda(J^{[q q_0]} / \aaa^{[q]} J^{[q q_0]})}{q^d}
= \lim_{q \rightarrow \infty} \frac{\lambda(K^{[q q_0]} / \aaa^{[q]} K^{[q q_0]})}{q^d}.
\]

On the other hand, by~\cite[Theorem 3.5(a)]{Vr}, we have \[
\lambda\left(\frac{K^{[q q_0]}} {\aaa^{[q]} K^{[q q_0]}}\right) 
 = \mu(K) \cdot \lambda\left( \frac{R}{\aaa^{[q]}} \right),
\]
and $\mu(K) = \ell^*(J)$.  These two equations displayed above,
then, give the result in case $J$ is $\m$-primary.
The fact that (\ref{eq:*spreadformula}) implies (\ref{eq:*spreadHKformula}) 
in this case is just by definition of Hilbert-Kunz multiplicities.

Now we drop the assumption that $J$ is $\m$-primary.  Let $\bfx = x_1, \dotsc, x_n$ be 
$R$-regular elements of $R$
whose images form a system of parameters for $R / J$.  By Lemma~\ref{lem:*indparam}, we can pick
an integer $t$ such that $K' := K + (\bx{t})$ is a $*$-independent ideal.
Moreover, $J' := J + (\bx{t}) \subseteq K^* + (\bx{t}) \subseteq (K + (\bx{t}))^* = {K'}^*$, so $K'$ is a minimal
$*$-reduction of $J'$, both of which are, of course, $\m$-primary.  What remains is to connect $J$ and $K$
with $J'$ and $K'$, respectively.

For each $i$ with $0 \leq i \leq n$, let $I_i = J + (x_1^t, \dotsc, x_i^t)$, and for $i <n$, $z_i = x_{i+1}^t$.
Then applying Lemma~\ref{lem:lengthinduction} to each $I = I_i$ and $z = z_i$ with $i<n$, we have \begin{align*}
\lim_{q \rightarrow \infty} \frac{\lambda(I_{i+1}^{[q q_0]} / \aaa^{[q]} I_{i+1}^{[q q_0]})}{q^d}
 &=  \lim_{q \rightarrow \infty} \frac{\lambda((I_i,z_i)^{[q q_0]} / \aaa^{[q]} (I_i,z_i)^{[q q_0]})}{q^d} \\
 &= \HKmult(\aaa) + \lim_{q \rightarrow \infty} \frac{\lambda(I_i^{[q q_0]} / \aaa^{[q]} I_i^{[q q_0]})}{q^d},
\end{align*}
so that, since $J' = I_n$ and $J = I_0$, after dividing by $\HKmult(\aaa)$ we have: \begin{align*}
\frac{1}{\HKmult(\aaa)} \lim_{q \rightarrow \infty} \frac{\lambda({J'}^{[q q_0]} / \aaa^{[q]} {J'}^{[q q_0]})}{q^d} 
 &= n + \frac{1}{\HKmult(\aaa)} 
 \lim_{q \rightarrow \infty}\frac{\lambda\left( J^{[q q_0]} / \aaa^{[q]} J^{[q q_0]} \right) } {q^d}.
\end{align*}

However, since $J'$ is $\m$-primary, we already know that the left hand side equals 
$\ell^*(J') = \mu(K') = n + \mu(K) = n + \ell^*(J)$.  Then subtracting $n$ from each side
gives the desired result.
\end{proof}

\section{Hilbert-Kunz multiplicity}

\begin{prop}\label{prop:HKrat}
Let $(R,\m)$ be an excellent analytically irreducible local ring such that $k = k(\bar{R})$, where
$\bar{R}$ is the normalization of $R$.  Let $I$ and $J$ be $\m$-primary ideals of $R$.  Then there
is some power $q_0$ of $p$ such that the following conditions are equivalent: \begin{enumerate}[\quad (a)]
\item\label{ratiq} There exist powers $q$, $q'$ of $p$ such that $q' \geq q \geq q_0$, and
$\HKmult(I J^{[q']})$ and $\HKmult(I J^{[q]})$ are both rational.
\item\label{ratbq} $\HKmult(I J^{[q]})$ is rational for all $q \geq q_0$.
\item\label{rat2} $\HKmult(I)$ and $\HKmult(J)$ are both rational.
\end{enumerate}
Moreover, there is some power $q_1$ of $p$ such that \[
\HKmult(J J^{[q]}) = (\ell^*(J) + q^d) \HKmult(J)
\]
for all $q \geq q_1$, where $d = \dim R$.  In particular, $\HKmult(J)$ is rational if and only if one such $\HKmult(J J^{[q]})$
is rational if and only if all such $\HKmult(J J^{[q]})$'s are rational.
\end{prop}

\begin{proof}
By Theorem~\ref{thm:*spreadlength}, there exists some $q_0$ such that for all $q \geq q_0$, we have
\begin{equation}\label{eq:a}
\ell \HKmult(I) + q^d \HKmult(J) = \HKmult(I J^{[q]}),
\end{equation}
where $\ell = \ell^*(J)$.  Hence, if $q' \geq q$ is another power of $p$, then we have
\begin{equation}\label{eq:b}
\ell \HKmult(I) + {q'}^d \HKmult(J) = \HKmult(I J^{[q']}),
\end{equation}
so that subtracting Equation~(\ref{eq:a}) from Equation~(\ref{eq:b}), we get: \begin{equation}\label{eq:J}
({q'}^d - q^d) \HKmult(J) = \HKmult(I J^{[q']}) - \HKmult(I J^{[q]}).
\end{equation}
On the other hand, if we multiply (\ref{eq:a}), by ${q'}^d$ and (\ref{eq:b}) by $q^d$, and then
subtract, we get: \begin{equation}\label{eq:I}
({q'}^d - q^d) \ell \HKmult(I) = {q'}^d \HKmult(I J^{[q]}) - q^d \HKmult(I J^{[q']}).
\end{equation}

It is trivial that $(\ref{ratbq}) \Rightarrow (\ref{ratiq})$.  Equation~(\ref{eq:a})
shows that $(\ref{rat2}) \Rightarrow (\ref{ratbq})$.  Equations~(\ref{eq:J})
and (\ref{eq:I}) show that $(\ref{ratiq}) \Rightarrow (\ref{rat2})$.

The second statement comes from replacing $I$ by $J$ in Equation~(\ref{eq:a}).
\end{proof}

The next Proposition does not refer to $*$-spread, but it is a nice base change formula for Hilbert-Kunz multiplicities
that works in a very general situation.
\begin{prop}\label{prop:HKBC}
Let $(R, \m) \rightarrow (S,\n)$ be a flat local homomorphism of Noetherian local
rings of prime characteristic $p>0$, such that $S / \m S$ is Cohen-Macaulay. Then
for any $\m$-primary ideal $\aaa$ in $R$ and any sequence $\z = z_1, \dotsc, z_s$ of
elements in $S$ whose images form a system of parameters for $S/\m S$, the following
two formulas hold: \begin{enumerate}
\item[(a)] \[
\lambda_S(S / (\aaa S, \z)) = \lambda_S (S / (\m S, \z)) \lambda_R (R / \aaa)
\]
\item[(b)] \[
\HKmul{S}(\aaa S + (\z)) = \mathrm{e}^{S / \m S}(\z) \HKmul{R}(\aaa).
\]
\end{enumerate}
\end{prop}

\begin{proof}
For part (a),
 we have that \[
S / (\aaa S, \z) \cong \frac{S / \z}{\aaa (S / \z)} \cong S/\z \otimes_R R / \aaa,
\]
and since $S / \z$ is flat over $R$, 
\begin{align*}
\lambda_{S / \z} (S / \z \otimes_R R / \aaa) 
 &= \lambda_{S / \z} ((S / \z) / \m(S / \z)) \cdot \lambda_R (R / \aaa) \\
 &= \lambda_{S / \m S} ((S / \m S) / \z(S / \m S)) \cdot \lambda_R (R / \aaa).
\end{align*}

For part (b), we replace $\z$ by $\z^{[q]}$ and $\aaa$ by $\aaa^{[q]}$ in (a) so that,
letting $d = \dim R$, we have \begin{align*}
\HKmul{S}(\aaa S + (\z)) &= \lim_{q \rightarrow \infty} \frac{\lambda_S(S / (\aaa S, \z)^{[q]})}{q^{d+s}} \\
&= \lim_{q \rightarrow \infty} \frac{\lambda_{S / \m S} ((S / \m S) / \z^{[q]}(S / \m S))}{q^s}
 \cdot \frac{\lambda_R (R / \aaa^{[q]})}{q^d} \\
&= \HKmul{S / \m S}(\z) \HKmul{R}(\aaa) = \mathrm{e}^{S / \m S}(\z) \HKmul{R}(\aaa).
\end{align*}
The last equality follows from \cite[Theorem 2]{Lech}.
\end{proof}

\section{$*$-spread and flat base change}

\begin{prop}\label{prop:*spreadBC}
Let $(R,\m) \rightarrow (S,\n)$ be a flat local homomorphism of prime characteristic $p>0$ 
excellent analytically irreducible Noetherian local rings which share a test element $c$. \begin{enumerate}
\item[(a)] If $x_1, \dotsc, x_n \in R$ are $*$-independent elements of $R$, they are
$*$-independent in $S$ as well.
\item[(b)] If $I$ is a proper ideal of $R$,
then $\ell^*(I) = \ell^*(I S)$.
\end{enumerate}
\end{prop}

\begin{proof}
For part (a), suppose that $x_n \in ((x_1, \dotsc, x_{n-1}) S)^*$.  Then 
for all $q \geq q_0$, $c x_n^q \in (x_1^q, \dotsc, x_{n-1}^q) S$, so that \begin{align*}
c &\in \bigcap_{q \geq q_0} (x_1^q, \dotsc, x_{n-1}^q) S :_S x_n^q
= \bigcap_{q \geq q_0} ((x_1^q, \dotsc, x_{n-1}^q) :_R x_n^q) S) \\
&= \left( \bigcap_{q \geq q_0} (x_1^q, \dotsc, x_{n-1}^q) :_R x_n^q \right) S,
\end{align*}
where both of the equalities follow from flatness of $S$ over $R$.  Now, 
$\bigcap_{q \geq q_0} ((x_1^q, \dotsc, x_{n-1}^q) :_R x_n^q)$ is contained in 
the union of the minimal primes of $R$, since $x_n \notin (x_1, \dotsc, x_{n-1})^*$.
That is, for some minimal prime $\p$ of $R$, $c \in \p S \cap R = \p$ by another
application of flatness, which is a contradiction.

As for part (b), let $J$ be a minimal $*$-reduction of $I$.  Let $x_1, \dotsc, x_n$ be
a minimal set of generators for $J$.  Then since $I S \subseteq J^* S \subseteq (J S)^*$ by
persistence of tight closure, $I S$ has a $*$-reduction generated by $\ell$ elements, which
shows that $\ell^*(I S) \leq \ell^*(I)$.\footnote{This part of the proof demonstrates a general fact
having nothing to do with any properties of the rings or the map between them.}  On the other hand,
$x_1, \dotsc, x_n$ are $*$-independent elements of $S$ by part (a), so $J S$ is a $*$-independent
ideal, so $\ell^*(I S) \geq \ell^*(I)$.
\end{proof}

\begin{cor}
Let $(R,\m,k) \rightarrow (S,\n,l)$ be a flat local homomorphism of Noetherian local rings $R$ and
$S$ which are excellent, analytically irreducible, and characteristic $p>0$.  Suppose further that
$R$ and $S$ share a test element $c$, that $k = \kappa(\bar{R})$ and $l = \kappa(\bar{S})$, and 
that $S / \m S$ is Cohen-Macaulay.

Then for any proper ideal $I$ of $R$ and any sequence $\z = z_1, \dotsc, z_r$ of elements of $S$
whose images in $S / \m S$ 
form a system of parameters for $S/\m S$, there is some power $q_0$ of $p$ such that \[
\lim_{q \rightarrow \infty} \frac{\lambda_S \left(
 (\m^{[q]}I^{[q q_0]} S + (\z)^{[q]} \cap I^{[q q_0]} S) / (\m S, \z)^{[q]} I^{[q q_0]} S \right)}{q^d} = 0,
\]
where $d = \dim S$.
\end{cor}

\begin{proof}
By Proposition~\ref{prop:HKBC}, $\HKmul{S}(\m S + \z) = \mathrm{e}^{S / \m S}(\z) \HKmult(\m).$

Also, we have \[
S / \z^{[q]} \otimes_R I^{[q q_0]} / \m^{[q]} I^{[q q_0]}
 \cong I^{[q q_0]} S / (\m^{[q]} I^{[q q_0]} S + (\z)^{[q]} \cap I^{[q q_0]} S),
\]
so that 
\[
\lambda_S(S / (\m S, \z^{[q]})) \lambda_R (I^{[q q_0]} / \m^{[q]} I^{[q q_0]})
= \lambda_S\left(\frac{I^{[q q_0]} S} {\m^{[q]} I^{[q q_0]} S + (\z)^{[q]} \cap I^{[q q_0]} S}\right)
\]
Letting $s = \dim S/\m S$, so that $\dim R = d - s$, we have: \begin{align*}
\lim_{q \rightarrow \infty}
 \lambda_S \left(\frac{I^{[q q_0]} S} {(\m S, \z)^{[q]} I^{[q q_0]} S} \right) /q^d
&= \HKmul{S}(\m S + \z) \ell^*(I S) \\
&= \mathrm{e}^{S/ \m S} (\z) \HKmul{R}(\m) \ell^*(I) \\
&= \lim_{q \rightarrow \infty} \frac{\lambda_S(S / (\m S, \z^{[q]}))}{q^s} \cdot
    \frac{\lambda_R (I^{[q q_0]} / \m^{[q]} I^{[q q_0]})}{q^{d-s}} \\
&= \lim_{q \rightarrow \infty} \lambda_S
 \left(\frac{I^{[q q_0]} S}{\m^{[q]} I^{[q q_0]} S + (\z)^{[q]} \cap I^{[q q_0]} S} \right) / q^d.
\end{align*}
Subtracting the right-hand side from the left-hand side gives the result.
\end{proof}


\begin{thebibliography}{99}
\bibitem[Ab]{Ab} I. Aberbach, {\em Extensions of weakly and strongly F-regular rings by flat maps}, J. Algebra, {\bf 241} (2001), 799--807.

\bibitem[Ep]{Epst} N. Epstein, {\em A tight closure analogue of analytic spread},  Math. Proc. Cambridge Philos. Soc., {\bf 139} (2005), no.2, 371--383.

\bibitem[HH]{HH} M. Hochster and C. Huneke, {\em Tight closure, invariant theory, and the Brian\c con-Skoda theorem},  J. Amer. Math. Soc., {\bf 3} (1990), 31--116.

\bibitem[Le]{Lech} C. Lech, {\em On the associativity formula for multiplicities}, Ark. Mat., {\bf 3} (1957), 301--314.

\bibitem[Mo]{Mo} P. Monsky, {\em The Hilbert-Kunz function}, Math. Ann. {\bf 263} (1983), 43--49.

\bibitem[NR]{NR} D. G. Northcott and D. Rees, {\em Reductions of ideals in local rings}, Proc. Cambridge Philos. Soc. {\bf 50} (1954), no. 2, 145--158.

\bibitem[Vr]{Vr}  A. Vraciu, {\em *-independence and special tight closure}, J. Algebra, {\bf 249} (2002), 544--565.
\end{thebibliography}
\end{document}